\documentclass[letterpaper, 10pt, conference]{ieeeconf}      
\IEEEoverridecommandlockouts                              
\usepackage{amsfonts,amsmath,amstext,mathrsfs}
\newcommand{\un}{\mathbf 1}
\newtheorem{theorem}{Theorem}[section]

\newtheorem{proposition}[theorem]{Proposition}
\newtheorem{corollary}[theorem]{Corollary}
\newtheorem{definition}[theorem]{Definition}

\renewcommand{\enspace}{\,}
\newcommand{\sP}{\mathcal{P}}
\newcommand{\NEW}[1]{{\em #1}\index{#1}}
\newcommand{\mrm}[1]{\text{#1}}
       
\def\<#1,#2>{\langle #1, #2\rangle}

\def\intt#1{\mathop{\mathrm{int}} (#1)}

\def\domain{\mathop{\mathrm{dom}}}
\def\lsc{\mathrm{lsc}}
\def\usc{\mathrm{usc}}
\def\supp{\mathop{\mathrm{idom}}}

\def\expf#1{\exp\left( #1\right)}
\def\quasi{\mathscr{QL}}
\def\maxlin{\mathscr{L}}
\def\devi#1{\mathscr{LD}_{#1}(Y)}
\def\dom#1{\domain (#1)}
\def\idom#1{\supp (#1)}

\def\ldom#1{\mathop{\mathrm{ldom}}(#1)}
\def\udom#1{\mathop{\mathrm{udom}}(#1)}
\def\sC{{\mathscr C}}
\def\sF{{\mathscr F}}
\def\sFC{{\mathscr F}_{\!c}}
\def\sM{{\mathscr M}}
\def\sG{{\mathscr G}}

\def\cP{{\mathscr P}}

\def\pW{\cP(W)}

\def\R{\mathbb{R}}
\def\N{\mathbb{N}}

\newcommand{\gal}{^{\circ}}
\newcommand{\bgal}{b\gal}
\newcommand{\set}[2]{\{#1\mid\,#2\}}

\newcommand{\rmaxb}{\overline{\R}_{\max}}

\newcommand{\rbar}{\overline{\R}}
\newcommand{\rmax}{\text{$\R_{\max}$}}%

\newcommand{\RBY}{\rbar^{_{\scriptstyle Y}}}

\newcommand{\RBX}{\rbar^{_{\scriptstyle X}}}

\newcommand{\Bd}{\mathop{\mathrm{bd}}(Y)}
\newcommand{\Cont}{\sC(Y)}
\newcommand{\Sci}{\lsc(Y)}
\newcommand{\Scs}{\usc(Y)}
\newcommand{\Contb}{\sC_{\mathrm{b}}(Y)}

\newcommand{\Scsb}{\usc_{\mathrm{b}}(Y)}
\newcommand{\For}{F}

\newcounter{myenumerate}
\renewcommand{\themyenumerate}{(\roman{myenumerate})}
\newenvironment{myenumerate}{\begin{list}{\themyenumerate }{
\usecounter{myenumerate}\setlength{\labelsep}{0.5ex}
\setlength{\leftmargin}{0pt}\setlength{\labelwidth}{-\labelsep}
}}{\end{list}}

\title{\LARGE \bf 
Solutions of max-plus linear equations and large deviations
\thanks{{\em Date}: September 13, 2005. Prepared for CDC-ECC'05.}}
\author{Marianne Akian, St\'ephane Gaubert and Vassili Kolokoltsov%
\thanks{Marianne Akian and St\'ephane Gaubert are with 
INRIA, Domaine de Voluceau,
B.P.~105, 78153 Le Chesnay Cedex, France.
{\tt  Marianne.Akian@inria.fr, Stephane.Gaubert@inria.fr}}
\thanks{Vassili Kolokoltsov is with 
Institute of the Information Transmission Problems of the 
Russian Academy of Sciences, Moscow, Russia.
{\tt vkolok@fsmail.net}}
}
 
\begin{document}
\maketitle
\addtolength{\textheight}{0.2cm}

\begin{abstract}
We generalise the G\"artner-Ellis theorem of large deviations theory.
Our results allow us to derive large deviation type results in stochastic 
optimal control from the convergence  of generalised 
logarithmic moment generating functions.
They rely on the characterisation
of the uniqueness of the solutions of max-plus linear equations.
We give an illustration for a simple investment model, in which 
logarithmic moment generating functions represent risk-sensitive values.
\end{abstract}

\section{Introduction}
Let $X,Y$ be two sets and $\rbar=\R\cup\{\pm\infty\}$
denote the extended real line. A Moreau conjugacy~\cite{moreau70},
associated to a kernel $b: X\times Y \to \rbar$,
is a map $B: \sF\to\sG$, where
$\sF$ and $\sG$ are subsets of $\RBY$ and $\RBX$ respectively, 
such that:
\[
Bf(x)= \sup\set{b(x,y)-f(y)}{y\in Y}, \quad \forall x\in X\enspace.\]
Here $b(x,y)-f(y)$ is an abbreviation of $b(x,y)+(-f(y))$,
with the convention that $-\infty$ is absorbing for addition.
An example of Moreau conjugacy
is the Legendre-Fenchel transform.
Moreau conjugacies are instrumental
in nonconvex duality, see~\cite[Chapter~11, Section~E]{rockwets},\cite{singer}.
The set $\rbar$ can be equipped with the semiring
structure of $\rmaxb$, in which the addition
is $(a,b)\mapsto \max(a,b)$ and the multiplication is
$(a,b)\mapsto a+b$, with the same convention as above.
If $B:\RBY\to\RBX$
is a Moreau conjugacy, the map $f\mapsto B(-f)$
is a max-plus linear operator with kernel.
Max-plus linear operators with kernel 
arise in deterministic optimal control and asymptotics,
and have been widely studied, see in 
particular~\cite{cuning,maslov92,BCOQ,maslovkolokoltsov95,DENSITE,gondran}.  

Given a map $g\in \sG$ and 
a Moreau conjugacy $B:\sF\to\sG$,
let us consider the problem:
\[
(\sP):\quad
\mrm{Find }f\in \sF\mrm{ such that } Bf =g\enspace,
\]
and more generally:
\[
(\sP'):\quad\begin{array}{l}
 \mrm{Find }f\in\sF\mrm{ such that }
Bf\leq g\\
\mrm{and }  Bf(x)=g(x)\mrm{ for all } x\in X'\enspace,
\end{array}
\]
where $X'\subset X$ is given.
In~\cite{AGK2} we gave effective conditions on $g$ for 
the solution $f$ of $(\sP')$ to exist and be unique, using 
generalised subdifferentials (associated to the Moreau conjugacy $B$).
We characterised the existence and uniqueness of the solution of $(\sP')$
in terms of coverings and minimal 
coverings by sets which are inverses of subdifferentials of $g$
(we recall these results in Section~\ref{moreau-sec}).
These conditions extend, to the case of infinite sets $X$ and $Y$,
the characterisations of the solutions of $(\sP)$
in terms of ``minimal resolvent coverings'' of $X$ which
were first shown by Vorobyev~\cite[Theorem  2.6]{vorobyev67} 
and then developed by Zimmermann~\cite[Chapter 3]{Zimmermann.K},
when $X$ and $Y$ are finite.

When $B$ is the Legendre-Fenchel transform, these results show
that essentially smooth convex functions have a
unique pre-image by the Legendre-Fenchel transform
(see~\cite[Corollary~6.4]{AGK2}), a fact which is
the essence of the classical G\"artner-Ellis theorem,
see e.g.~\cite[Theorem~2.3.6,(c)]{DEMBO} for a general presentation.
Indeed, as we shall see in Section~\ref{sec-quasiconv}, 
Problem $(\sP')$ arises naturally
when looking for the rate function in large deviations.

Large deviation type asymptotics arise in optimal control when studying
the long-term behaviour of some controlled stochastic process.
For instance, assume that some real nonnegative controlled diffusion $X_t$,
representing the wealth of an investor, has an asymptotic growth rate,
which means that
 $\lim_{t\to +\infty} \frac{\log(X_t)}{t}$ exists almost surely,
and assume that this growth rate satisfies a large deviation principle 
with factor $1/t$, which means in loose terms that for
``good'' subsets $A$ of $\R$, $-\frac{1}{t} 
\log P( \frac{\log(X_t)}{t}\in A)$ tends
to $I(A):=\inf_{x\in A} I(x)$ for some rate function $I:\R\to [0,+\infty)$.
Then, one may want to find a control minimising the rate $I(A)$.
This problem was studied by Pham in~\cite{pham}, to which the reader is 
referred for more information.

In the present paper, we show how the results of~\cite{AGK2}
yield a characterisation of the rate function of a large deviation
principle, generalising the G\"artner-Ellis theorem
(Section~\ref{sec-quasiconv}).
We also study rate functions in optimal stochastic control,
such as the maximal long term growth rate of wealth described above.
To develop these results in a unified way, we introduce
 (in Section~\ref{sec-quasi}) the notion of quasi (max-plus) linear form.
It includes as special cases, possibly up to
a log-exp conjugacy: integrals with respect to finite measures,
suprema of such integrals and idempotent integrals
with respect to finite cost measures~\cite{DENSITE}.
We then introduce (in Section~\ref{sec-weak})
the notion of weak convergence of quasi-linear forms,
which generalises the large deviation principle of Varadhan.
Finally, we illustrate the results on a simple investment model
(Section~\ref{sec-appli}).

\section{Images and pre-images of Moreau conjugacies}
\label{moreau-sec}
We briefly recall some of the results of~\cite{AGK2}.

\subsection{Moreau conjugacies}
Let $X$ and $Y$ be two Hausdorff topological spaces.
Denote by $\sF$ the set of lower semicontinuous (l.s.c.)
maps from $Y$ to $\rbar$ and by $\sG$ the set $\RBX$ of
all maps from $X$ to $\rbar$.
The sets $\sF$ and $\sG$ are endowed with the partial ordering of
functions, for which they are complete lattices.
Let $b: X\times Y\to \rbar$ be a map which is lower semicontinuous
in the second variable. 
Then the maps $B:\sF\to \sG$ and
$B\gal:\sG\to\sF$ defined by
\begin{align*}
Bf(x)&= \sup \set{b(x ,y)-f(y)}{y\in Y}
\; \forall f\in \sF, \; x\in X\enspace,\\ 
B\gal g(y)&= \sup \set{b(x,y)-g(x)}{x\in X} 
\; \forall g\in \sG,\; y\in Y 
\end{align*}
are Moreau conjugacies~\cite{moreau70} 
and they are conjugate to each other, or in duality, meaning that 
$(B,B\gal)$ is a dual Galois connection
(see~\cite[Theorem~2.1 and Example~2.7]{AGK2}).
Moreover, by~\cite[Theorem~2.1]{AGK2}, the map $b$
is uniquely determined by the Moreau conjugacy $B$.
We call it the \NEW{kernel}
of $B$. The kernel of $B\gal$ is necessarily equal to
the symmetric map of $b$, denoted by $\bgal$:
$Y\times X\to \rbar,\; (y,x)\mapsto \bgal(y,x)=(x,y)$.
Taking two topological vector spaces $X$ and $Y$ in duality,
and $b(x,y)=\<x,y>$, we obtain the classical Legendre-Fenchel 
transform $Bf=f^*$.

In the sequel we shall assume that 
$b(x,y)\in\R\cup\{-\infty\}$ for all $(x,y)\in X\times Y$, 
and that for all $x\in X$, there exists $y\in Y$ such that
$b(x,y)\in\R$, and symmetrically that for all $y\in Y$, there
exists $x\in X$ such that $b(x,y)\in \R$.

\subsection{Existence of pre-images}

We shall use the following notion of subdifferentials of Moreau conjugacies
(see~\cite{balder,dolecki,lindberg,martinez88,martinez95}).
Given  $f\in\sF$ and $y\in Y$, 
the \NEW{subdifferential} of $f$ at $y$ with respect to $b$ (or $B$),
denoted by $\partial^b f(y)$, or $\partial f(y)$ for brevity, is
the set:
\begin{eqnarray*}
\lefteqn{\partial f(y)= \{ x\in X\mid b(x,y)\in\R,}\\
&&b(x,y')-f(y') \leq b(x,y)-f(y)\; \forall y'\in Y\}\enspace.
\end{eqnarray*}
For $g\in \sG$ and $x\in X$, the subdifferential of $g$ at $x$
with respect to $\bgal$,  $\partial^{\bgal} g(x)$, will be denoted 
by $\partial\gal g (x)$ for brevity.
When $b(x,y) =\<x,y>$
we recover the classical definition of subdifferentials. 

We shall use the following covering notions.
When $\Phi$ is a map from a set $Z$ to the set $\pW$ of all subsets of some
set $W$, we set 
$\Phi^{-1}(w)=\set{z\in Z}{w\in \Phi(z)}$.
If $Z'\subset Z$ and $W'\subset W$,
$\{\Phi(z)\}_{z \in Z'}$ is a {\em covering}
of $W'$ if $\cup_{z\in Z'} \Phi(z)\supset W'$. 
An element $y\in Z'$ is {\em algebraically essential} with respect to
this covering 
if $W'\not\subset \cup_{z\in Z'\setminus\{y\}} \Phi(z)$.
When $Z$ is a topological space,
$y$ is {\em topologically essential} 
if, for all open neighbourhoods
$U$ of $y$ in $Z'$, $W'\not\subset \cup_{z\in Z'\setminus U} \Phi(z)$.
The covering of $W'$ by $\{\Phi(z)\}_{z \in Z'}$ is 
{\em algebraically (resp.\ topologically) minimal} if all elements
of $Z'$ are algebraically (resp.\ topologically) essential.


The kernel $b$ is said \NEW{coercive} 
if for all $x\in X$, and all neighbourhoods $V$ of $x$ in $X$, the function 
\begin{align}\label{defi-coerc}
y\in Y\mapsto b_{x,V}(y)=
\sup_{z\in V} b(z,y)-b(x,y) 
\end{align}
has relatively compact finite sublevel sets, which
means that $\set{y\in Y}{b_{x,V}(y)\leq \beta}$
is relatively compact for all $\beta\in \R$. 
We also denote by $\sFC$ the set of all $f\in \sF$ such that
for all $x\in X$, $y\mapsto b(x,y)-f(y)$ has
relatively compact finite superlevel sets, which
means that for all $\beta\in\R$, the set 
$\set{y\in Y}{b(x,y)-f(y)\geq \beta}$ is relatively compact.
For any map $g$ from a topological space $Z$ to $\rbar$,
we set: 
$ \ldom g := \set{z\in Z}{g(z)<+\infty}$,
$\udom g := \set{z\in Z}{g(z)>-\infty}$,
$\dom g  := \ldom g \cap \udom g$ (the \NEW{domain} of $g$),
$\idom g =  \set{z\in \dom g}{
\limsup_{z'\to z} g(z')<+\infty}$.

We shall occasionally need the following assumptions:
\catcode`\@=11
\def\refcounter#1{\protected@edef\@currentlabel
       {\csname the#1\endcsname}%
}
\catcode`\@=12
\newcounter{assume}
\def\theassume{{\rm (A\arabic{assume})}}
\def\theassumep{{\rm (A\arabic{assume})$'$}}
\def\myitem{\refstepcounter{assume}\item[\theassume\hskip 0.72ex]}
\def\myitemp{\refcounter{assumep}\item[\theassumep]}

\begin{itemize}[\setlabelwidth{(A4)'}]
\myitem\label{a-yisdiscrete} $Y$ is discrete;

\myitemp\label{a-bgalgisfinite} $b$ is continuous in the second variable,
and $B\gal g(y)>-\infty$ for all $y\in Y$; 

\myitem\label{a-bgalginfc}\label{casei} $B\gal g\in\sFC$;

\myitemp\label{caseii}
\label{a-biscoercive}
$b$ is coercive and $X'\subset \idom g\cup g^{-1}(-\infty)$.

\myitem\label{globalassump}
Conditions~\ref{a-yisdiscrete} or~\ref{a-bgalgisfinite},
and~\ref{a-bgalginfc} or~\ref{a-biscoercive} hold.
\end{itemize}

\begin{theorem}[\protect{\cite[Theorem 3.5]{AGK2}}]\label{cover}
Let $X'\subset X$ and $g\in \sG$.
Consider the following statements:
\begin{myenumerate}
\item\label{e-eq1}Problem $(\sP')$ has a solution,
\item \label{e-c2} $\{(\partial\gal g)^{-1}(y)\}_{y\in \ldom{B\gal g}}$
is a covering of $X' \cap \udom g$.
\end{myenumerate}
We have 
\ref{e-c2}$\Rightarrow$\ref{e-eq1}. 
If~\ref{globalassump} is satisfied, 
then \ref{e-eq1}$\Leftrightarrow$\ref{e-c2}.
\end{theorem}

\subsection{Uniqueness of the pre-image}

A map $h$ from a topological space $Z$ to $\rbar$ 
is \NEW{quasi-continuous}~\cite{NEUBRUNN}
if for all open sets $G$ of $\rbar$, the set $h^{-1}(G)$ 
is included in the closure of its interior.
When $h$ is l.s.c., this is equivalent to the condition
that $h$ is the l.s.c.\ hull of the upper semicontinuous (u.s.c.) hull of $h$.

\begin{theorem}[\protect{\cite[Theorem 4.6]{AGK2}}]\label{minicover1}
Let $X'\subset X$ and $g\in \sG$.
Assume that $\{ (\partial\gal g)^{-1}(y)\}_{y\in \ldom{B\gal g}}$ 
is a covering of $X' \cap \udom g$, and denote by
$Z_a$ (resp.\ $Z_t$) the set of algebraically (resp.\ topologically)
essential elements with respect to this covering.
Let $Z=Z_a\cup \intt{Z_t}$, where $\intt{Z_t}$ denotes the interior of $Z_t$
relatively to $\dom{B\gal g}$.
Assume that~\ref{globalassump} is satisfied
and that $B\gal g$ is quasi-continuous on its domain.
Then Problem  $(\sP')$ has a solution, and 
any solution  $f$ of $(\sP')$ satisfies
\[
f\geq B\gal g,\quad\mrm{and}\quad f(y)=B\gal g(y) \quad
\mrm{ for all }y\in Z\enspace.
\]
\end{theorem}
\begin{theorem}[\protect{\cite[Theorem 4.7]{AGK2}}]\label{minicover}
Let $X'\subset X$ and $g\in \sG$.
Consider the following statements:
\begin{myenumerate}
\item\label{fe-eq1} Problem $(\sP')$ has a unique solution,
\item\label{fe-c1} $\{ (\partial\gal g)^{-1}(y)\}_{y\in \ldom{B\gal g}} $
is a topologically minimal covering of $X' \cap \udom g$.
\end{myenumerate}
If~\ref{globalassump} is satisfied, then \ref{fe-eq1}$\Rightarrow$\ref{fe-c1}.
If in addition $B\gal g$ is quasi-continuous on its domain,
then \ref{fe-eq1}$\Leftrightarrow$\ref{fe-c1}.
\end{theorem}

\section{Quasi-linear forms}\label{sec-quasi}

We assume now that $Y$ is a Polish (complete separable metric) space.
We denote by $\rmax$ the subsemiring of $\rmaxb$ composed of the
elements of $\R\cup\{-\infty\}$.
The set $(\rmax)^Y$ is a sublattice of $\RBY$, which is 
conditionally complete.
We denote by $\sup$ or $\vee$ (resp.\ $\inf$ or  $\wedge$)
the supremum (resp.\ infimum) operation.
The set $(\rmax)^Y$ can be endowed with the semimodule structure in
which the addition is $(f,g)\mapsto f\vee g$,
and the scalar multiplication is $(a,f)\in \rmax\times (\rmax)^Y
\mapsto a+f\in (\rmax)^Y$ with $(a+f)(y)=a+ f(y)$ for all $y\in Y$.
(Semimodules and subsemimodules are defined as modules and 
submodules over rings~\cite{litvinov00,cgq02}).
We denote by $\Bd$ (resp.\ $\Cont$, resp.\ $\Sci$, resp.\ $\Scs$)
the set of functions from $Y$ to
$\rmax$  that are bounded above by a constant 
(resp.\ continuous, resp.\ l.s.c., resp.\ u.s.c.).
We also use the notations
$\Contb:=\Cont\cap\Bd$ and $\Scsb:=\Scs\cap\Bd$.
All these sets are subsemimodules of $(\rmax)^Y$.

\begin{definition}
Let $\sM$ be a subsemimodule of $(\rmax)^Y$.
A map $\For:\sM\to \rmax$ (or $\rmaxb$) is a
\NEW{quasi- (max-plus) -linear form} (on $\sM$) if it is isotone, that is
\begin{subequations}\label{quasilinear}
\begin{equation}\label{quasilinear.0}
 \varphi\leq \psi\implies \For(\varphi)\leq \For(\psi)\;\text{for all }\varphi,\psi\in \sM\enspace,
\end{equation}
if it is additively homogeneous, that is
\begin{equation}
 \For(\lambda +\varphi)= \lambda + \For(\varphi)\;\text{for all }
\lambda\in\rmax,\; \varphi\in \sM\enspace,\label{quasilinear.1}
\end{equation}
and if there exists $\alpha\in \rmax$ such that
\begin{equation}
\For(\varphi\vee \psi)\leq \alpha + \For(\varphi)\vee\For(\psi)
\;\text{for all }\varphi,\psi
\in \sM \enspace .\label{quasilinear.2}
\end{equation}
A quasi-linear form $\For$ on $\sM$ is \NEW{continuous}
if it preserves nondecreasing  converging sequences.
\end{subequations}
\end{definition}
We denote by $\rho(\For)$ the infimum of the $\alpha$
satisfying~\eqref{quasilinear.2}, and 
by $\quasi(\sM)$ the set of continuous quasi-linear forms from $\sM$
to $\rmax$.
When $\For$ takes at least one value in $\R$, $\rho(\For)\geq 0$ and 
one can take $\alpha=\rho(\For)$ in~\eqref{quasilinear.2}. 
Otherwise $\rho(\For)=-\infty$.


A map $\For:\sM\to\rmaxb$ is a quasi-linear form such that 
$\rho(\For)\leq 0$ if, 
and only if, $\For$ is a max-plus linear form, that is $\For$ 
satisfies~\eqref{quasilinear.1} and 
$\For(\varphi\vee \psi)=
 \For(\varphi)\vee\For(\psi)$ 
for all $\varphi,\psi \in \sM$. 
We set $\maxlin (\sM):=\set{\For\in\quasi (\sM)}{\rho(\For)\leq 0}$.
Given $f:Y\to \rbar$, 
the map $\For: (\rmax)^Y\to\rmaxb$ 
defined by
\begin{equation}\label{max-dens}
\For (\varphi)=\sup_{y\in Y} \left(\varphi(y)-f(y)\right) 
\text{ for all } \varphi\in (\rmax)^Y 
\end{equation}
is a continuous (as a quasi-linear form)
max-plus linear form on $(\rmax)^Y$ 
and so on $\Contb$. 
A map $f$ satisfying~\eqref{max-dens}
is called a \NEW{density} of $\For$.
Conversely, since $Y$ is a separable metric space,
any element of $\maxlin(\Contb)$ 
has a unique l.s.c.\ density~\cite[Th. 4.8 and Cor. 3.12]{DENSITE} (see
also~\cite{KOLO88,kolokoltsov92,maslovkolokoltsov95}),
which is bounded below by some real constant.
Note however that a max-plus linear form is not necessarily continuous.

Let $\mu$ be a finite positive measure on $Y$, let $\varepsilon>0$
and consider the map $\For:\Contb\to\rmax$ with 
\begin{equation}\label{large-devi}
\For (\varphi)=
\varepsilon \log \left( \int_Y \expf{\frac{\varphi(y)}{\varepsilon}}d\mu(y)
\right)\enspace , 
\end{equation}
for all $\varphi\in \Contb$.
Then $\For$ is a continuous quasi-linear form with $\rho(\For)\leq
 \varepsilon\log (2)$.
The maps $F$ of the form~\eqref{large-devi} where $\mu$ is
a probability measure occur in large deviations principles.
We shall denote by $\devi{\varepsilon}$ the set of all such maps.

Let $\For_i\in\quasi(\sM)$, for $i\in I$, such that 
$\sup_{i\in I} \rho(\For_i)<+\infty$. Then the map 
$\sup_{i\in I} \For_i:\sM\to\rmaxb,\; \varphi\mapsto
 \sup_{i\in I} \For_i (\varphi)$
is a continuous quasi-linear form on $\sM$ 
and it satisfies $\rho(\sup_{i\in I} \For_i)\leq \sup_{i\in I} \rho(\For_i)$.

\begin{proposition}\label{extens}
Any $\For\in\quasi(\Contb)$ 
admits a unique extension to a continuous quasi-linear
form on $\Sci$ (with values in $\rmaxb$) that we also denote by $\For$:
\[
 \For(\varphi)=\sup_{\psi\in\Contb,\; \psi\leq \varphi} \For(\psi)
\quad \text{ for all } \varphi\in\Sci \enspace ,
\]
and a maximal extension to a continuous quasi-linear
form on $(\rmax)^Y$ that we also denote by $\For$:
\[
 \For(\varphi)=\inf_{\psi\in\Sci,\; \psi\geq \varphi} \For(\psi)
\quad \text{ for all } \varphi\in (\rmax)^Y  \enspace .
\]
The value of $\rho(\For)$ for the maximal extension of $\For$ to $(\rmax)^Y$
and for its restriction to $\Contb$ coincide.
\end{proposition}

If $A$ is  a subset of $Y$, we denote by  $\un_A:Y\to\rmax$ 
the max-plus characteristic 
function of $A$: $\un_A(y)=0$ if $y\in A$ and $\un_A(y)=-\infty$ 
otherwise.
If $\For$ is as in Proposition~\ref{extens}, we shall also denote by 
$\For$ the map 
\( 
\For :\sP(Y)\to\rmax,\quad A\mapsto \For (\un_A)\).
This map is isotone: $A\subset B\implies \For (A)\leq \For (B)$,
and it 
satisfies some inner and outer-continuity properties.
If $\For$ is a continuous max-plus linear form 
or is an element of $\devi{\varepsilon}$,
with $\varepsilon>0$, then the exponential of its restriction to
$\sP(Y)$ is a capacity in the sense defined in~\cite{VERV91,VERV95}.
Some other related sets of functions on $\sP(Y)$ are defined
in~\cite{jiang,charalamb}.

\section{Weak convergence of quasi-linear forms}\label{sec-weak}
Let $\For\in \quasi(\Contb)$ and
$(\For_n)_{n\in \N}$ be a sequence of $\quasi(\Contb)$. 
We say that $\For_n$ \NEW{weakly converges} towards $\For$ if 
\( \lim_{n\to\infty}
\For_n(\varphi)=\For(\varphi)\mrm{ for all }\varphi\in\Contb\). 
In that case, we get that $\rho(\For)\leq \liminf_{n\to\infty} \rho(\For_n)$.
When $\For_n \in \devi{\varepsilon}$ is defined from the  measure $\mu_n$,
then $\For\in \devi{\varepsilon}$  and the weak convergence
of $\For_n$ towards $\For$  is equivalent to the weak (or narrow) convergence 
of $\mu_n$ towards the measure $\mu$ corresponding to $\For$.
If all the $\For_n$ are continuous max-plus linear forms
 with l.s.c.\ densities $f_n$, then $\For$ is also a continuous
max-plus linear form and the weak convergence
of $\For_n$ towards $\For$ is equivalent to the weak convergence of the
cost measure with density $f_n$ towards the cost measure with density
$f$ (see~\cite{DUALITY}), where $f$ is the l.s.c.\ density of $\For$.

We say that an element $\For$ of $\quasi(\Contb)$
is \NEW{tight} if \( \inf_{K\subset Y,\; K\text{ compact}} \For
(K^{c})=-\infty\). 
A sequence $(\For_n)_{n\in \N}$ of $\quasi(\Contb)$
is \NEW{asymptotically tight}
if $\limsup_{n\to \infty}\rho(\For_n)<+\infty$ and
\[ \inf_{K\subset Y,\; K\text{ compact}} \limsup_{n\to\infty}
\For_n(K^{c})=-\infty\enspace .\]
Since $Y$ is a Polish space, any element of $\devi{\varepsilon}$,
with $\varepsilon>0$, is tight. 
An element of $\maxlin(\Contb)$ is tight if,
and only if, its l.s.c.\ density $f$ is inf-compact, that is
$\set{y\in Y}{f(y)\leq \alpha}$ is compact for all $\alpha\in\R$.


\begin{theorem}\label{defi-equiv}
Let $\For\in \quasi(\Contb)$ and
$(\For_n)_{n\in \N}$ be a sequence of $\quasi(\Contb)$. 
Denote also by $\For$ and $\For_n$ the extensions given by
Proposition~\ref{extens}.
Consider the following statements:
\begin{align}
\label{defi-equiv.1}& \For_n \text{ weakly converges towards } \For\enspace,\\
\label{defi-equiv.2}& \liminf_{n\to\infty}\For_n(\varphi)\geq \For(\varphi)  \text{ for all }
\varphi\in \Sci \enspace,\\
\label{defi-equiv.3}& \limsup_{n\to\infty}\For_n(\varphi)\leq \For(\varphi)  \text{ for all }
\varphi\in \Scsb \enspace,\\
\label{defi-equiv.4}& \liminf_{n\to\infty}\For_n(G)\geq \For(G)  \text{ for all open } G \subset Y\enspace,\\
\label{defi-equiv.5}& \limsup_{n\to\infty}\For_n(C)\leq \For(C) 
 \text{ for all closed } C \subset Y\enspace, \\
\label{defi-equiv.6}& \limsup_{n\to\infty}\For_n(K)\leq \For(K) 
 \text{ for all compact } K \subset Y  \enspace.
\end{align}
We have {\rm (\ref{defi-equiv.2},\ref{defi-equiv.3})}%
$\Rightarrow$\eqref{defi-equiv.1}$\Rightarrow$\eqref{defi-equiv.2}%
$\Rightarrow$\eqref{defi-equiv.4},
\eqref{defi-equiv.3}$\Rightarrow$\eqref{defi-equiv.5}$\Rightarrow$%
\eqref{defi-equiv.6}.
If $\rho(\For)\leq \lim_{n\to\infty} \rho(\For_n )=0$, then~%
\eqref{defi-equiv.1}$\Leftrightarrow$%
{\rm (\ref{defi-equiv.2},\ref{defi-equiv.3})}$\Leftrightarrow$%
{\rm (\ref{defi-equiv.4},\ref{defi-equiv.5})}.
If $(\For_n)_{n\in\N}$ is asymptotically
tight then~\eqref{defi-equiv.5}$\Leftrightarrow$\eqref{defi-equiv.6}.
\end{theorem}

When $\For_n\in\devi{\varepsilon_n}$ with corresponding 
probability measure $\mu_n$,
$\For$ is a continuous max-plus linear form with  density $f$,
and $\lim_{n\to\infty} \varepsilon_n=0$,
$(\mu_n)_{n\in\N}$ obeys the \NEW{large deviation principle}
of Varadhan~\cite{VAR} with rate function $f$ if, and only if,
$f$ is nonnegative and 
inf-compact and~{\rm (\ref{defi-equiv.4},\ref{defi-equiv.5})} holds.
In this context, the
implication~{\rm (\ref{defi-equiv.4},\ref{defi-equiv.5})}$\Rightarrow$\eqref{defi-equiv.1} is called the contraction
principle of Varadhan, and some other implications in
Theorem~\ref{defi-equiv} are proved in~\cite{VAR} and 
in~\cite[Theorem 3.1.3]{tolya01} (see also~\cite{DEMBO} and~\cite{BRYC90}).
In the context of capacities, the conditions
{\rm (\ref{defi-equiv.4},\ref{defi-equiv.6})} define the vague
convergence and  the conditions
{\rm (\ref{defi-equiv.4},\ref{defi-equiv.5})}  define the 
narrow (weak) convergence~\cite{VERV91}.
In the context of continuous max-plus linear forms, some of the
implications in Theorem~\ref{defi-equiv} are proved in~\cite{DUALITY}.
The following result is also classical for large deviations.
It was stated for max-plus linear forms in~\cite{DUALITY}.

\begin{theorem}\label{compact-tight}
Let $(\For_n)_{n\in \N}$ be a sequence of $\quasi(\Contb)$ 
such that $\limsup_{n\to\infty}\For_n(Y)<+\infty$ and
$\lim_{n\to\infty} \rho(\For_n )=0$.
Then  there exists  
$\For\in\maxlin(\Contb)$ and a subsequence
of $(\For_n)_{n\in \N}$ such that 
{\rm (\ref{defi-equiv.4},\ref{defi-equiv.6})} holds for that subsequence.
\end{theorem}

\section{Uniqueness of pre-images of Moreau conjugacies and convergence 
of quasi-linear forms}\label{sec-quasiconv}

Let $X,Y,B,B\gal,b,b\gal$ be as in Section~\ref{moreau-sec}.

We say that $b$ is \NEW{strongly coercive} if 
for all $x\in X$ and all neighbourhoods $V$ of $x$ in $X$,
there exists a finite subset $W$ of $V$ such that
the function $b_{x,W}$ defined as in~\eqref{defi-coerc}
has relatively compact finite sublevel sets.

We say that $b$ is \NEW{upper (strongly) coercive} if 
for all $x\in X$, and all neighbourhoods $V$ of $x$ in $X$,
there exists a finite subset $W$ of $V$ such that
$b(x,\cdot)$ is bounded above on each finite sublevel set of $b_{x,W}$.

If $b$ is strongly coercive and continuous in the second variable,
 then $b$ is coercive and upper coercive.
If $b(x,y)=\<x,y>$, then $b$ is upper coercive
(take $W=\{t x\}$ with $t>1$ near enough from $1$).
If in addition $X=Y=\R^n$, then $b$ is strongly coercive
(take $W=\set{x\pm \varepsilon e_i}{1\leq i\leq n}$
with $\varepsilon>0$ small enough, and $(e_1,\ldots, e_n)$ a
basis of $\R^n$).

The following result motivates the study of Problem $(\sP')$.
\begin{theorem}\label{th-ineq}
Let $(\For_n)_{n\in \N}$ be a sequence of $\quasi(\Contb)$ 
such that $\lim_{n\to\infty} \rho(\For_n )=0$. 
 Assume that $\For_n$ weakly converges towards
$\For\in\maxlin(\Contb)$, with
l.s.c.\ density $f$, and that $b$ is continuous in the second variable
and upper coercive. Let $g:X\to\rbar$ be defined by:
\begin{equation}\label{th-ineq-as}
g(x)=\limsup_{n\to \infty} \For_n(b(x,\cdot))\quad \text{ for all } 
x\in X\enspace .
\end{equation}
Then
\begin{equation}\label{pp}
 Bf\leq g \text{ and } Bf=g \text{ on } \idom g\cup g^{-1}(-\infty)
\enspace .\end{equation}
\end{theorem}

The next result follows from Theorems~\ref{defi-equiv},
\ref{compact-tight}, \ref{th-ineq}, \ref{cover}, \ref{minicover1}
and~\ref{minicover}. It needs the following technical assumption:
\begin{itemize}[\setlabelwidth{(A4)}]
\myitem\label{technicalassump}
Conditions~\ref{a-yisdiscrete} or~\ref{a-bgalgisfinite} hold;
Condition~\ref{a-bgalginfc} holds or $b$ is coercive; $b$ is upper coercive;
and $B\gal g$ is quasi-continuous on its domain.
\end{itemize}

\begin{theorem}\label{coro-gartner}
Let $(\For_n)_{n\in \N}$ be an asymptotically tight sequence 
of $\quasi(\Contb)$, such that $\limsup_{n\to\infty}\For_n(Y)<+\infty$
and $\lim_{n\to\infty} \rho(\For_n )=0$.
Let $g:X\to \rbar$ be defined by~\eqref{th-ineq-as} and denote
by $\overline{\For}$ the continuous max-plus linear form
with density $B\gal g$. 
Assume that~\ref{technicalassump} is satisfied. Then
\begin{myenumerate}
\item \label{coro-gartner.1} There exists 
$\For\in\maxlin(\Contb)$, and a subsequence of $(\For_n)_{n\in \N}$ 
which converges weakly towards $\For$.
\item \label{coro-gartner.2}
$\{ (\partial\gal g)^{-1}(y)\}_{y\in \ldom{B \gal g}}$
is a covering of $\idom g$.
\item \label{coro-gartner.3}
If $\For$ is an accumulation point of  $(\For_n)_{n\in \N}$ 
for the weak convergence, and if $f$ is the l.s.c.\ density of $\For$,
then $f\geq B\gal g$. Hence
$\limsup_{n\to\infty}\For_n(C)\leq \overline{\For}(C)$
for all closed $C\subset Y$.
\end{myenumerate}
Assume in addition that the limsup in~\eqref{th-ineq-as} is a limit,
and let $Z$ be defined as in Theorem~\ref{minicover1}
with $X'=\idom g$.

\begin{myenumerate}\addtocounter{myenumerate}{3}
\item \label{coro-gartner.4}
If $\For$ and $f$ are as in~\ref{coro-gartner.3},
then $f=B\gal g$ on $Z$. Hence 
$\liminf_{n\to\infty}\For_n(G)\geq \overline{\For}(G\cap Z)$ for all open
$G\subset Y$.
\item \label{coro-gartner.5}
 If $\{ (\partial\gal g)^{-1}(y)\}_{y\in \ldom{B \gal g}}$
is a topologically minimal covering of $\idom g$,
then $\For_n$ weakly converges towards $\overline{\For}$.
\end{myenumerate}
\end{theorem}

The following result can be used to obtain the ``compactness'' of the
sequence $(\For_n)_{n\in \N}$.

\begin{proposition}\label{tighness}
Let $(\For_n)_{n\in \N}$ be a sequence of $\quasi(\Contb)$ 
such that $\lim_{n\to\infty} \rho(\For_n )=0$, and 
let $g:X\to\rbar$ be given by~\eqref{th-ineq-as}.
Assume that $b$ is strongly coercive, and
that there exists $x_0\in\idom g$ such that $b(x_0,\cdot)$ is
bounded below by some real constant.
Then $\limsup_{n\to\infty} \For_n(Y)<+\infty$ and 
$(\For_n)_{n\in \N}$ is asymptotically tight.
\end{proposition} 

\begin{corollary}[Generalised G\"artner-Ellis theorem]
\label{coro-gartnerbis}
Let $(\For_n)_{n\in \N}$ be a sequence of $\quasi(\Contb)$ 
such that $\lim_{n\to\infty} \rho(\For_n )=0$, 
let $g:X\to\rbar$ be given by~\eqref{th-ineq-as}, 
and assume that the limsup there is a limit.
Assume that~\ref{a-yisdiscrete} or \ref{a-bgalgisfinite}
hold, that $b$ is strongly coercive, that 
$B\gal g$ is quasi-continuous on its domain and
that there exists $x_0\in\idom g$ such that $b(x_0,\cdot)$ is
lower bounded by some real constant.
Then the conclusions of Theorem~\ref{coro-gartner} hold.
\end{corollary} 

When $B$ is the Legendre-Fenchel transform on $\R^n$, $x_0=0$, $\For_n\in
\devi{\varepsilon_n}$ with $\lim_{n\to\infty} \varepsilon_n=0$,
the statement of Corollary~\ref{coro-gartnerbis} contains 
the G\"artner-Ellis theorem as stated in~\cite[Th.~2.3.6]{DEMBO}.
Indeed, $b(0,\cdot)\equiv 0$, $b$ is strongly coercive
(see above),
$B\gal g$ is quasi-continuous on its domain (see~\cite[Lemma 6.1]{AGK2}),
$B\gal g(y)>-\infty$ for all $y\in Y$ when
$g$ is  proper, in particular when $\idom g\neq \emptyset$.
Moreover, by~\cite[Proposition 6.3 and Corollary~6.4]{AGK2},
$\{ (\partial\gal g)^{-1}(y)\}_{y\in \ldom{B \gal g}}$
is a topologically minimal covering of $\idom g$, when
$g$ is an essential smooth  l.s.c.\ proper convex function on $\R^n$, 
which means that the interior of its domain $\idom g$ is nonempty, that
$g$ is differentiable in  $\idom g$, and that
the norm of the differential of $g$ at $x$ tends  to infinity when
$x$ goes to the boundary of $\dom g$, see~\cite[Section~26]{ROCK}.

The proof of our generalisation of 
the G\"artner-Ellis theorem essentially uses compactness
arguments together with the uniqueness
of the pre-image of an essential smooth convex function
by the Legendre-Fenchel transform.
This last argument was made explicit by O'Brien 
and Vervaat~\cite[Theorem~4.1 (c)]{VERV95},
and Puhalskii~\cite[Lemmas 3.2 and 3.5]{PUHAL94-1}
for the G\"artner-Ellis theorem, and by
Gulinsky~\cite[Theorems 4.7 and 5.3]{gulinski}
for the more general case where $B$ is the Legendre-Fenchel transform, and 
$\For_n(\varphi)=\varepsilon_n \log J_n(\exp (\frac{\varphi}{\varepsilon_n}))$ 
with $J_n(\varphi\vee\psi)\leq J_n(\varphi)+J_n(\psi)$.

If $\For_n$ weakly converges towards $\For$ with a density
$f$ which is not essentially strictly convex
(or equivalently such that its Legendre-Fenchel transform
$f^*$ is not essentially smooth)
the G\"artner-Ellis theorem only gives
the inequalities of Assertions~\ref{coro-gartner.3}
and~\ref{coro-gartner.4} 
of Theorem~\ref{coro-gartner} with $Z\neq Y$, thus
the rate function $f$ cannot be identified.
The classical method is to adapt the proof of the G\"artner-Ellis theorem,
whereas using Theorem~\ref{coro-gartner}, one may simply consider 
a different kernel $b$ than that of the Legendre-Fenchel transform.
Moreover, Proposition~\ref{tighness} can be applied to another kernel.

The following result is useful in the study of optimal control problems.
\begin{theorem}\label{sup-gartner}
Assume that, for $i\in I$, $(\For_{n,i})_{n\in \N}$ is an
asymptotically tight sequence of $\quasi(\Contb)$, such that 
$\limsup_{n\to\infty} \For_{n,i}(Y)<\infty$ 
and $\lim_{n\to\infty} \rho(\For_{n,i})=0$.
Let $g:X\to\rbar$ be defined by
\[
g(x)=\sup_{i\in I} \limsup_{n\to \infty} \For_{n,i}(b(x,\cdot))
\quad \text{ for all } x\in X\enspace ,
\]
and denote
by $\overline{\For}$ the continuous max-plus linear form
with density $B\gal g$. 
Assume that~\ref{technicalassump} is satisfied. Then
\begin{myenumerate}
\item \label{sup-gartner.1} There exists 
$\For\in\maxlin(\Contb)$ such that 
$\sup_{i\in I} \limsup_{n\to \infty} \For_{n,i}(G)\geq \For(G)$
 for all open $G\subset Y$, and $\sup_{i\in I} \limsup_{n\to \infty}
 \For_{n,i}(C)\leq \For(C)$ for all closed $C\subset Y$.
\item \label{sup-gartner.2}
The l.s.c.\ density $f$ of $\For$ satisfies~\eqref{pp}.
Hence $\{ (\partial\gal g)^{-1}(y)\}_{y\in \ldom{B \gal g}}$
is a covering of $\idom g$.
\item \label{sup-gartner.3}
We have $f\geq B\gal g$, hence, for all closed $C\subset Y$,
$\sup_{i\in I} \limsup_{n\to \infty} \For_{n,i}(C)\leq  \overline{\For}(C)$.
\item \label{sup-gartner.4}
Let $Z$ be defined as in Theorem~\ref{minicover1} for $X'=\idom g$.
Then  $f=B\gal g$ on $Z$. Hence, for all open $G\subset Y$,
$\sup_{i\in I} \limsup_{n\to \infty} \For_{n,i}(G)\geq 
\overline{\For}(G\cap Z)$. 
\end{myenumerate}
\end{theorem}

\section{An application to the optimal long-term rate of an
investment model}
\label{sec-appli}
We consider here the simple Merton model~\cite{merton}
of an investor who has the possibility to invest in one bank account
paying a fixed  interest rate $r>0$
and in one stock or risky asset whose price is a log-normal diffusion
with  expected rate $\alpha>r$ and rate variation $\sigma$,
and who has the ability to transfer funds between the assets with no cost.
We denote by $W_t$ the total wealth of the investor at time $t$,
and by $\xi_t$ the proportion of fund invested in the risky asset.
The process $W_t$ satisfies the following stochastic differential equation:
\[
dW_t= (r+(\alpha-r)\xi_t) W_t dt +\sigma \xi_t W_t dB_t\enspace ,
\]
where $B_t$ is a Brownian motion.
The control process $\xi=(\xi_t)_{t\geq 0}$ is supposed to be
adapted to the Brownian filtration and stationary.
We allow borrowing and shortselling, which means that 
$\xi_t$ can be any real number.
One is interested in
maximising some function of the long term growth rate 
of the investor.
One possibility is to consider the risk-sensitive problem
\[
\sup_{\xi} \limsup_{T\to\infty} 
\frac{1}{T(1-\gamma)} \log E[(W_T)^{1-\gamma}]
\]
where $E$ denotes the expectation and 
$\gamma$ is the risk-aversion coefficient.
Another possibility is to consider, for $c\in\R$:
\begin{equation}\label{seuil}
\sup_{\xi} \limsup_{T\to\infty} 
\frac{1}{T} \log P[(\log (W_T)/T\geq c ]\enspace .\end{equation}
In~\cite{pham} the latter 
problem was considered for a different investment model
and the relation with the risk-sensitive problems with $\gamma<1$ 
was discussed and used to obtain a result of the same nature as
the G\"artner-Ellis theorem.

We apply here the results of the previous sections
to compute the quantity~\eqref{seuil}.
Let $Y=X=\R$ and consider the quasi-linear form $\For_{W_0,T,\xi}$
on $\Contb$ defined by
\[ \For_{W_0,T,\xi}(\varphi)=\frac{1}{T} \log 
E[\exp(T \varphi(\log(W_T)/T))\mid W_0]
\enspace \]
and extended as in Proposition~\ref{extens},
together with the quasi-linear form $\For_{W_0,T}=\sup_{\xi} \For_{W_0,T,\xi}$.
Then,  for all $W_0,T,\xi$, $\For_{W_0,T,\xi}(Y)=\For_{W_0,T}(Y)=0$, and
$\For_{W_0,T,\xi}$ and $\For_{W_0,T}\in\quasi(\Contb)$.
Moreover $\rho(\For_{W_0,T,\xi})\leq \log(2)/T$, 
thus $\lim_{T\to\infty} \rho(\For_{W_0,T,\xi})=
\lim_{T\to\infty} \rho(\For_{W_0,T})=0$.

Let $b(x,y)=xy$ be the kernel of the Legendre-Fenchel transform.
Then for all $x\in\R$, $\For_{W_0,T,\xi}(b(x,\cdot))=
\frac{1}{T} \log E[(W_T)^x\mid W_0]$ is a risk-sensitive 
utility function.
We have the homogeneity property: $\For_{W_0,T,\xi}(b(x,\cdot))=
\frac{x\log(W_0)}{T}+\For_{1,T,\xi}(b(x,\cdot))$.
Let $g: X\to\rbar$ be defined by 
\begin{equation}\label{defiginvest}
g(x)= \sup_{\xi\in\R} x\left(r+ (\alpha-r) \xi +
(x-1)\frac{\sigma^2 \xi^2}{2}\right)
\end{equation}
for $x\in X$. Then 
$g(x)=x ( r+\frac{(\alpha-r)^2}{2 \sigma^2 (1-x)}) $
if $0\leq x<1$ and $g(x)=+\infty$ otherwise.
Moreover, for $0\leq x<1$, the proportion 
\( \bar{\xi}_x= \frac{\alpha-r}{\sigma^2 (1-x)}\)
realises the maximum in~\eqref{defiginvest}.
We have, for all $T>0$ and $x\in X$,
$\For_{1,T}(b(x,\cdot))=g(x)$,
and for  all $T>0$, $W_0>0$  and $0\leq x<1$, 
the constant control process $\xi_t\equiv  \bar{\xi}_x$
maximises $\For_{W_0,T,\xi}(b(x,\cdot))$.
Hence
\begin{eqnarray}
\lefteqn{\sup_{\xi}\limsup_{T\to\infty} \For_{W_0,T,\xi}(b(x,\cdot))}
\label{suplimsup}\\
&&= \lim_{T\to\infty}\For_{W_0,T}(b(x,\cdot))=g(x) \enspace.\nonumber 
\end{eqnarray} 
The Legendre-Fenchel transform $g^*$ of g is given by 
$g^*(y)=(\sqrt{y-r}-\frac{\alpha-r}{\sqrt{2}\sigma})^2$ if 
$y\geq z_0:=r+\frac{(\alpha-r)^2}{2 \sigma^2}$
and $g^*(y)=0$ otherwise. So,
if  $\overline{\For}$ and $Z$ are defined 
as in Theorem~\ref{coro-gartner}, we get $Z=(z_0,+\infty)$, and
$\overline{\For}((c,+\infty)\cap Z)=\overline{\For}([c,+\infty))=
-g^*(c)$ for all $c\in\R$.
If, for any sequence $T_n$ going to infinity,
the sequence $(\For_{W_0,T_n})_{n\geq 0}$
were asymptotically tight, then Theorem~\ref{coro-gartner}
would show:
\begin{subequations}\label{investconclutheo}
\begin{align}
&\!\liminf_{T\to +\infty} \For_{W_0,T}(G)\geq \overline{\For}(G\cap Z) 
\;\text{for all open}\; G \subset \R\enspace,\\
&\!\limsup_{T\to +\infty} \For_{W_0,T}(C)\leq \overline{\For}(C)
\;\text{for all closed}\; C \subset \R\enspace.
\end{align}
\end{subequations}
In particular this would show:
\begin{eqnarray}
\lefteqn{\lim_{T\to +\infty} \sup_{\xi} \frac{1}{T}
\log P[(\log (W_T)/T\geq c \mid W_0]}\nonumber\\
&&= \lim_{T\to +\infty} \For_{W_0,T}([c,+\infty))=-g^*(c)\enspace.
\label{mertonconclu}
\end{eqnarray}
However, since $0\not\in\idom g$, one cannot use 
Proposition~\ref{tighness} to show the asymptotic tightness 
of $(\For_{W_0,T_n})_{n\geq 0}$.

Let us thus replace the process $\log(W_T)/T$
by its maximum with some constant $a<z_0$. This amounts to
replacing $\For_{W_0,T,\xi}$ by the quasi-linear form 
$G_{W_0,T,\xi}(\varphi) :=\For_{W_0,T,\xi}(\varphi\circ \chi_a)$
where $\chi_a(x)=x\vee a$. We also consider
$G_{W_0,T}=\sup_{\xi} G_{W_0,T,\xi}$.
We take now $Y=[a,+\infty)$ and $X=[0,+\infty)$.
The kernel $b(x,y)=xy$ is strongly coercive with respect to these
new sets $X$ and $Y$. The corresponding Moreau conjugacies $B$ and $B\gal$
are the Legendre-Fenchel transform composed with the restriction operation
to $Y$ and $X$ respectively.
Since $\rho(\For_{W_0,T,\xi})$ tends to $0$ when $T$ goes to infinity,
and $b(x,\chi_a(y))=b(x,y)\vee xa$ for all $x\in [0,+\infty)$
 and $y\in \R$, we get that 
$\lim_{T\to\infty}G_{W_0,T}(b(x,\cdot))=
\lim_{T\to\infty}\For_{W_0,T}(b(x,\cdot))\vee xa=g(x)$
for all $x\in X$.
Moreover,  $B\gal g$ is the restriction of $g^*$ to $Y$.
With respect to the new set $Y$,
$b(x,\cdot)$ is lower bounded for all $x\in X$ and
$\idom g=[0,1)$, hence Proposition~\ref{tighness} shows that 
$G_{W_0,T_n}$ is asymptotically tight for any sequence $T_n$ tending
to infinity. Then the conclusions~\eqref{investconclutheo}
of Theorem~\ref{coro-gartner} hold with $\For_{W_0,T}$
replaced by $G_{W_0,T}$, and with $Z=(z_0,+\infty)$ unchanged.
Since $\un_A\circ \chi_a=\un_A$ if $a\not\in A$, and $a$ can be chosen
small enough, we deduce~\eqref{mertonconclu}.

Let us now apply Theorem~\ref{sup-gartner} to
the sequences $(G_{W_0,T_n,\xi})_{n\in \N}$, with $T_n$ tending to infinity,
and $W_0$ fixed, and where the parameter $i$ corresponds to 
the couple composed of the control process
$\xi$ and of the sequence $(T_n)_{n\geq 0}$.
We obtain, by the same arguments as before, that for all $c\in\R$:
\begin{eqnarray}
\lefteqn{\sup_{\xi} \limsup_{T\to +\infty} \frac{1}{T}
\log P[(\log (W_T)/T\geq c \mid W_0]}\nonumber\\
&&= \sup_{\xi} \limsup_{T\to +\infty} \For_{W_0,T,\xi}([c,+\infty))
=-g^*(c)\enspace.
\label{mertonsup}
\end{eqnarray}
The latter conclusion is of the same nature as the one 
of~\cite[Theorem 3.1]{pham}.
Note however that for the proof of~\eqref{mertonsup} one does not need 
that the supremum in~\eqref{suplimsup} is attained and that for this maximum
the limsup is a limit, as is required in~\cite{pham}, even 
if these properties hold in our example. However, these conditions
were useful to prove~\eqref{mertonconclu}, and so
to prove that in~\eqref{mertonsup} 
the sup and limsup operations commute.

\end{document}